\documentclass[11pt]{amsart}
\usepackage{amsmath, amscd}
\usepackage{amsfonts}
\usepackage{amssymb}
\usepackage{amsthm}
\usepackage{verbatim}
\usepackage[mathscr]{euscript}
\usepackage[toc,page]{appendix}
\usepackage{xspace}
\usepackage{fullpage}
\usepackage{graphicx}

\def\ul{\underline}
\def\ol{\overline}
\def\ds{\displaystyle}
\newcommand{\R}{\ensuremath{\mathbb{R}}\xspace}

\newcommand{\Pf}{{\bigskip\noindent\bf Proof.\quad} }

\newcommand{\wtilde}{\widetilde}

\newtheorem{theorem}{Theorem}[section] 
\newtheorem{lemma}[theorem]{Lemma} 
\newtheorem{proposition}[theorem]{Proposition}
\newtheorem{corollary}[theorem]{Corollary}
\theoremstyle{definition} 

\newtheorem{remark}[theorem]{Remark}

\begin{document}
\title{On Vorticity Gradient growth for the axisymmetric 3D euler equations without swirl}
\author{Tam Do}
\email{tamdo@usc.edu}
\address{Dept. of Mathematics, University of Southern California, Los Angeles, CA, USA 90089}
\date{\today}
\maketitle

\begin{abstract}
We consider solutions of the 3D axisymmetric Euler equations without swirl. In this setting, well-posedness is well-known due to the essentially 2D geometry. The quantity $\omega^\theta/r$ plays an analogous role as vorticity in 2D. For our first result, we prove that the gradient of $\omega^\theta/r$ can grow with at most double exponential rate with improving a priori bound close to the axis of symmetry. Next, on the unit ball, we show that at the boundary, one can achieve double exponential growth of the gradient of $\omega^\theta/r$.
\end{abstract}

\section{Introduction}

In this paper, we examine vorticity gradient growth in 3D {\it axisymmetric flows without swirl}. Global well-posedness of the 3D Euler equations for this class of flows is shown in the work of Ukhovsksii and Yudovich \cite{Yudovich} (see \cite{majda2002vorticity, Raymond, Shirota} for related results). This class of flows shares similarities with two-dimensional flows where well-posedness is also well-known \cite{majda2002vorticity, Holder, Marchioro, Wolibner}. In two dimensions, vorticity is conserved along particle trajectories while in 3D axisymmetric flows without swirl, the quantity $\omega^\theta(r,z)/r$ plays an analogous role where $\omega^\theta$ is the angular component of vorticity in cylindrical coordinates and $r$ being the radial variable. This fact is key in showing well-posedness for this class of 3D flows as this provides a priori bounds for vorticity.

\smallskip
Our first main result concerns an upper bound on the growth of the gradient of $\omega^\theta/r$.  For 2D flows, the upper bound for the gradient of vorticity is double exponential growth in time \cite{Yudovich1}. We will prove a similar upper bound in the axisymmetric case. However, we also show that this upper bound improves to essentially exponential growth near the axis. We make this precise in the next section. This result will serve in contrast to 2D Euler flows as we rule out any double exponential growth
of the gradient at the axis. The special structure of the axisymmetric Biot-Savart law is used.

\smallskip
For our second result, we explore the sharpness of this upper bound and construct an example of double exponential gradient growth. Such growth will occur on the boundary of the unit ball $B(0,1)=\{(r,z): \, r^2+z^2\le 1\}$ away from the axis. For the 2D Euler equations, there have been a number of recent results concerning the gradient growth of vorticity \cite{KiselevSverak, Denisov1, Denisov2, Zlatos}. The techniques we will use bear most resemblance to those of Kiselev and Sverak \cite{KiselevSverak} who construct an example of double exponential vorticity gradient growth on the boundary of a unit disk. Their initial data is inspired by
the ``singular cross" of Bahouri and Chemin \cite{BahouriChemin} and the authors show 
that particle trajectories are approximately hyperbolic near the desired point of gradient growth. We will construct an initial data inspired by this scenario for the ball. In order to prove that trajectories have
such structure, we require a closed-form expression for the Green's function of an elliptic operator
to be specified below. This is one of the primary reasons the ball is chosen for our domain rather than the more natural choice of a cylinder. We anticipate that our construction will work on other domains such as a cylinder.

\smallskip
For both of our results, the proofs rely on adequate expressions and estimates for the Biot Savart law for axisymmetric flows without swirl. In the axisymmetric setting, the Biot Savart law is considerably more complicated than the law for fluid velocity $u$ and vorticity $\omega$ in 2D Euler which is
\begin{align*}
u(x,t) = \nabla^\perp \int_D G_D(x,y)\omega(y,t) dy.
\end{align*}
Here, $G_D$ is Green's function for the Laplacian for the Dirchlet problem on a 2D simply connected domain $D$. In section 2.2, we make a precise statement of the Biot Savart law we use. Similar to 2D Euler, one can express the velocity in terms of the vorticity $\omega^\theta$ integrated against some kernel. Away from the axis, this kernel has similar estimates as $\nabla^\perp G_D$, but near the axis, the kernel will have better decay estimates. This similarity away the axis will lead to the double exponential growth at the boundary away from the axis. Additionally, the better kernel decay for points near the axis will lead to our improved upper bound referenced above.

\smallskip
We believe our work is a small step toward bridging ideas of small scale creation in 2D to 3D. For the 3D axisymmetric Euler equations {\it with} swirl, a potential scenario for singularity formation was proposed Luo and Hou \cite{HouLuo1} based upon their numerical simulations. A singularity is reported on the boundary of a cylinder and flow is observed to have hyperbolic structure. In our setting without swirl, we produce an example with double exponential gradient growth where the flow has hyperbolic-type structure. Proving singularity formation for Euler flows with swirl would require many deep new ideas. Another interesting question is the possible singularity formation of the 3D axisymmetric Euler equations with swirl at the axis of symmetry. 

\section{The setup}
Consider the 3D Euler equations for a velocity field $u$ and pressure $p$
\begin{align}
\label{euler1}
u_t+ u\cdot \nabla u +\nabla p  &=0 \\
\nabla\cdot u &=0 \\
u(\cdot, 0) &= u_0
\end{align}
on $D\times (0,\infty)$ where $D$ is either the unit ball $B(0,1)$ or a 
finite radius cylinder with periodic boundary condition in $z$.
In addition, we have no-flow condition on the solid boundary:
\begin{align*}
u\cdot n =0\quad \mbox{on} \quad \partial D.
\end{align*}
Here, we consider $u$ which is axisymmetric without swirl. Specifically, the velocity field $u$ will have the form
\begin{align*}
u(r,z,t)=u^r(r,z,t)\, {\bf e_r}+ u^z(r,z,t)\,{\bf  e_z}
\end{align*}
where ${\bf e_r}=\begin{bmatrix} \cos \theta & \sin \theta & 0 \end{bmatrix}^t$ and ${\bf e_z}=\begin{bmatrix} 0 & 0 & 1 \end{bmatrix}^t$ are unit vectors from cylindrical coordinates, the other unit vector being ${\bf e_\theta}= \begin{bmatrix} -\sin \theta & \cos\theta & 0\end{bmatrix}^t$. Since the vector field $u$ has no ${\bf e_\theta}$ component, the flow $u$ is said to be without swirl. After a change of coordinates, $\eqref{euler1}$ 
becomes
\begin{equation}
\label{euler}
(\partial_t + u^r \partial_r+ u^z \partial_z)w=0, \quad\quad w(r,z,t)= \frac{\omega^\theta(r,z,t)}{r}
\end{equation}
where $\omega^\theta$ is the angular component of vorticity $\omega :=\mbox{curl}\, u$. Throughout, we will assume $w_0(r,z)=w(r,z,0)\in L^\infty$ and so $\|w(\cdot,t)\|_{L^\infty}\le \|w_0\|_{L^\infty}$. The relations between vorticity, velocity, and stream function are as follows
\begin{align}
\label{vorticity}
\omega^\theta=u_z^r-u_r^z, \quad u^r = -\frac{\psi_z}{r}, \quad u^z = \frac{\psi_r}{r}
\end{align}

which can be combined to get
\begin{align}
\label{law}
L\psi := -\frac{\psi_{rr}}{r^2}+\frac{\psi_r}{r^3}-\frac{\psi_{zz}}{r^2}= \frac{\omega^\theta}{r}=w.
\end{align}
Often times, it will be convenient to use an equivalent form of \eqref{law} which is
\begin{align}
\label{law1}
\mathscr{L}\left(\frac{\psi}{r}\right) := -\left(\partial_r^2 +\frac{1}{r}\partial_r +\partial_z^2 - \frac{1}{r^2}\right) \frac{\psi}{r} =\omega^\theta.
\end{align} 
From these relations, in the next section, we will derive the Biot-savart law relating $\omega^\theta$ and $u$. Thus, $\omega^\theta$ completely describes the flow.

\medskip
Due to axial symmetry, there are additional conditions at the axis. To ensure solutions remain smooth, the stream function must satisfy
\begin{align}
\label{pole}
 \partial_r^{(2m)} \left.\left( \frac{\psi}{r}\right) \right|_{(0+,z)}=0 \quad\mbox{$m=0,1,2,\ldots$}
\end{align}
which implies the following conditions on $u$
\begin{align}
\partial_r^{(2m+1)} u^z(0+,z)=0, \quad \partial_r^{(2m)} u^r(0+,z)=0 \quad\mbox{$m=0,1,2,\ldots$}.
\end{align}
These conditions follow from the below lemma from Liu and Wang \cite{LiuWang}:

\begin{lemma}
\label{wang}
{\bf (a)}
Let $u$ be a $C^k$ smooth (in Cartesian coordinates) 3D axisymmetric vector field $u(r,\theta, z)=u^r(r,z) {\bf e_r} +u^\theta(r,z) {\bf e_\theta} + u^z(r,z){\bf e_z}$. Then $u^r, u^\theta, u^z\in C^k([0,\infty)\times \R)$ and 
\begin{align*}
\partial_r^{2\ell +1} u^z(0^+,z) = 0, \quad 1\le 2\ell + 1\le k, \\
\partial_r^{2m} u^r(0^+,z)=\partial_r^{2m} u^\theta(0^+,z) = 0, \quad 0\le 2m\le k.
\end{align*}

\bigskip
{\bf (b)} Suppose $\phi(r,z)\in C^{k+1}([0,\infty),\R), f\in C^k([0,\infty),\R)$ satisfying $\partial_r^{2m} \phi(0^+,z)=0$ for $0\le 2m\le k+1$ and 
$\partial_r^{2\ell} f(0^+,z)=0$ for $0\le 2\ell \le k$. Then the vector field 
\begin{align*}
u := -\partial_z \phi\, {\bf e_r}+ \frac{\partial_r(r\phi)}{r}\, {\bf e_z} + f\, {\bf e_\theta}
\end{align*}
for $r>0$ is a $C^k$ smooth (in Cartesian coordinates) 3D axisymmetric vector field with a removable singularity at $r=0$.

\end{lemma}

\noindent
The above lemma has a direct analogoue for vector fields over the domains we consider. Also, we note that the incompressibility condition becomes
\begin{align}
\label{divfree}
(ru^r)_r+ (ru^z)_z=0
\end{align}
in cylindrical coordinates.
Recall the similarity between the system \eqref{euler}, \eqref{vorticity}, and \eqref{law} with the 2D Euler equation in vorticity form:
\begin{align*}
\omega_t+ u\cdot \nabla \omega =0, \quad -\Delta\psi= \omega, \,\,\, u= \nabla^\perp \psi.
\end{align*}

\subsection{Statement of Main Results}

We will consider the system on the domain $D$ subject to no-flow boundary condition on solid boundaries
\begin{align*}
u\cdot n =0\quad\mbox{on}\quad \partial D.
\end{align*} 
Again, we will only consider the case when $D$ is the unit ball or a finite radius cylinder with periodic boundary condition in $z$. Much of our work will generalize to other axially symmetric domains. We indicate specifically in 
the proofs below where our specific domain choice is used.

\medskip
Define the trajectory map $\Phi_t(r,z)=(\Phi_t^r(r,z,t), \Phi_t^z(r,z,t))$ associated with \eqref{euler} by
\begin{align*}
\frac{d}{dt} \Phi_t(r,z) = u(\Phi_t(r,z),t), \quad \Phi_0(r,z)=(r,z).
\end{align*}
As the existence of classical solutions to \eqref{euler} is known \cite{Yudovich, majda2002vorticity, Raymond, Shirota},
the goal of our work is to address the sharpness of a priori bounds for such solutions.
Our first main result is the following upper bound on the gradient of $w$.

\begin{theorem}
\label{upper}
Let $\ds w_0= \frac{\omega^\theta(r,z,0)}{r}\in C^1(D)$ and $\omega_0^\theta(0,z)=0$ for $(0,z)\in \ol{D}$. 
Let $w(r,z,t)$ be the corresponding classical solution of \eqref{euler}.

\medskip
{\bf(a)}
 We have the following double exponential in time growth estimate for the gradient of $w(r,z,t)$ 
\begin{align}
\label{upper1}
1+ \log\left(1+ \frac{\|\nabla w(\cdot,t)\|_{L^\infty}}{\|w_0\|_{L^\infty}}\right) \le \left(1+  \log\left(1+ \frac{\|\nabla w_0\|_{L^\infty}}{\|w_0\|_{L^\infty}}\right)\right)\exp(C\|w_0\|_{L^\infty }t).
\end{align}
for some constant $C$.

\medskip
{\bf(b)} Particle trajectories can only approach the axis of symmetry with at most exponential rate. That is, there
exists a constant $C$ dependent on $D$ and $w_0$ such that for every $x\in D$
\begin{align}
\label{lowb}
|\Phi_t^r(x)| \ge r \exp(-Ct)
\end{align}

\medskip
{\bf (c)} The solution $w(r,z,t)$ satisfies the following estimate
\begin{align*}
\sup_{r>0} \frac{|w(r,z,t)-w(0,z,t)|}{r} \le \|\nabla w_0\|_{L^\infty} \exp(Ct)
\end{align*}

\end{theorem}

For the 2D Euler equations, estimates similar to $\eqref{lowb}$ have been shown when additional symmetry is
imposed on the vorticity or the domain has a corner (\cite{Ito1, Elgindi, Ito2}). Next, we provide an example of double exponential growth at the boundary.

\begin{theorem}
\label{lower}
Consider the 3D axisymmetric Euler equations without swirl on the unit ball $B(0,1)$. There exists initial data $w_0$ such that $\|\nabla w_0\|_{L^\infty(B(0,1))}/\|w_0\|_{L^\infty(B(0,1))}>1$ and the solution $w(r,z,t)$ of \eqref{euler} satisfies the following lower bound
\begin{align*}
\frac{\|\nabla w(r,z,t)\|_{L^\infty}}{\|w_0\|_{L^\infty}} \ge \left( \frac{\|\nabla w_0\|_{L^\infty}}{\|w_0\|_{L^\infty}}\right)^{C\exp (C\|w_0\|_\infty t)}
\end{align*}

\end{theorem}

\medskip

\subsection{Biot-Savart law}

\smallskip

Here, we will give a somewhat general derivation of the axisymmetric Biot-Savart law. Later on, we require rather explicit expressions for $u$ and its derivatives in terms of $w$. We'll emphasize the analogy with the 2D case rather than the general 3D Biot Savart law on domains which is highly non-trivial \cite{Enciso}. In order to derive the relation between $u$ and $\omega^\theta$, one must solve the following system
\begin{align}
\label{biotsavart} 
u_z^r-u_r^z &=\omega^\theta \\
\mbox{div}\, u & = 0 \\
u\cdot n &= 0 \quad\mbox{on} \,\, \partial D 
\end{align}
Recall, by the divergence-free condition \eqref{divfree}, there exists a scalar stream function $\psi(r,z)$ such that
 
\begin{align}
\label{psiperp}
 u^r= - \frac{\psi_z}{r}\quad \mbox{and}\quad  u^z= \frac{\psi_r}{r}
\end{align}

Now, let us find a vector potential $A$ such that
\begin{align}
\label{vectorid}
u(r,z,t)= u^r \, {\bf e_r} + u^z\, {\bf e_z} = \mbox{curl}\, A, \quad \mbox{div}\, A=0
\end{align}
Such an $A$ must satisfy
\begin{align}
\label{Aeq}
\omega^\theta\, {\bf e_\theta} = \mbox{curl}\, u = \mbox{curl}\, \mbox{curl}\, A = - \Delta A.
\end{align}

\medskip
Now, using the Dirchlet Green's function for the (scalar) Laplacian in 3D, we can write an expression for the Biot Savart law.
First, we write the Green's function for the Laplacian as 
\begin{align}
\label{3DGreen}
G_D(x,y)= \frac{1}{4\pi|x-y|}+h(x,y)
\end{align}
where for each $y\in D$ the corrector function $h$ solves
\begin{align*}
\Delta_x h(x,y) &= 0 \\
h(x,y) &=-\frac{1}{4\pi|x-y|}, \,\, x\in \partial D
\end{align*}

Writing $y$ in cylindrical coordinates $y=(r,\theta,z)$ and without loss of generality asumming $x=(\bar{r},0,\bar{z})$, we have
\begin{align*}
|x-y| = \sqrt{r^2-2r\bar{r}\cos \theta +\bar{r}^2+(z-\bar{z})^2}.
\end{align*}
Without loss of generality we can refer to $D$ as our given axially symmetric domain in $\R^3$ or express it in coordinates $(r,z)$, $r\ge 0$, $z\in \R$, depending on context.  This abuse of notation can be justified since we are in the axisymmetric setting so quantities depend only on their values on one $\theta$ plane of $D$. In cylindrical $(r,z)$ coordinates, $\partial D$ will only be points that correspond to boundary points in 3D. That is, the points on the axis that are not boundary points in 3D do not ``become" boundary points once in $(r,z)$ coordinates.
Then we write the ${\bf e_\theta}$ component of $A$ as
\begin{align*}
A^\theta(\bar{r},\bar{z})=\int_D \mathscr{A}^\theta(\bar{r},\bar{z},r,z) \omega^\theta(r,z) \, dr\, dz
\end{align*}
where we define
\begin{align}
\label{eq:A}
\mathscr{A}^\theta(\bar{r},\bar{z},r,z) = \frac{1}{4\pi}\int_0^{2\pi} r\cos\theta\left(\frac{1}{\sqrt{r^2-2r\bar{r}\cos\theta+\bar{r}^2+(z-\bar{z})^2}}+h(\bar{r},\bar{z},r,\theta,z)\right) \, d\theta.
\end{align}
The $r\cos\theta$ factor above comes from $\eqref{Aeq}$ where the vector Laplacian is used.
By the following lemma,  we have a relation between
$A$ and $\psi$.

\begin{lemma}
\label{lawform}
	Consider \eqref{Aeq} with boundary condition $A|_{\partial D} =0$. Suppose $\omega^\theta(r,z)\in C^1$
	and $\omega^\theta(0,z)=0$. Then the vector field $A$ only has ${\bf e_\theta}$ component. In addition, the ${\bf e_\theta}$ component of $A$ is $\theta$-independent. Moreover, $A$ has the following form
\begin{align}
\label{Apsi}
	A= \frac{\psi(r,z)}{r}\, e_\theta
\end{align}
where $\psi$ satisfies 
\begin{align}
\label{psi_bdy}
\psi&=0\quad \mbox{on}\,\, \partial D \\
\label{psi_eq}
\mathscr{L} \left(\frac{\psi}{r}\right)&:= -\left(\partial_r^2 +\frac{1}{r}\partial_r +\partial_z^2 - \frac{1}{r^2}\right) \frac{\psi}{r}=\omega^\theta.
\end{align}
\end{lemma}

\Pf 

By Lemma \ref{wang} and our hypothesis on $\omega^\theta$, $\omega^\theta\, {\bf e}_\theta$ is a continuous and bounded vector field and we have solvability of 
equation $\eqref{Aeq}$. 

\medskip
First, consider the case when $D$ is the unit ball $B(0,1)$. We can just directly compute.
Recall the Green's function for the ball in cylindrical coordinates
\begin{align}
\label{Green_ball}
G_B(r,z,\bar{r},\bar{z})  =\frac{1}{4\pi\sqrt{r^2-2r\bar{r}\cos\theta+\bar{r}^2+(z-\bar {z})^2}}- \frac{1}{4\pi\sqrt{r^2+z^2} \sqrt{(r^\ast)^2-2r^\ast\bar{r}\cos\theta+\bar{r}^2+(z-\bar{z})^2}}
\end{align}
where $\ds r^\ast = \frac{r}{r^2+z^2}$ and $\ds z^\ast = \frac{z}{r^2+z^2}$. We can use the Green's function $\eqref{Green_ball}$ to solve \eqref{Aeq}.
When integrated against $\sin \theta$, the contribution from the first term of the Green's function $\eqref{Green_ball}$ to $A$ is zero as
$$
\int_0^{2\pi} \frac{r\sin\theta}{\sqrt{r^2-2r\bar{r}\cos\theta+\bar{r}^2+(z-\bar{z})^2}}\, d\theta = 0.
$$
The other term of the Green's function can be handled similarly so 
\begin{align*}
A(\bar{r},\bar{z})  &= \int_D G_B(r,z,\bar{r},\bar{z}) \omega^\theta(r,z) \begin{bmatrix} -\sin \theta \\ \cos\theta \\ 0\end{bmatrix} \, r\, d\theta\, dr\, dz  \\
&=  \int_D G_B(r,z,\bar{r},\bar{z}) \omega^\theta(r,z) \begin{bmatrix}0 \\ \cos\theta \\ 0\end{bmatrix} \, r\, d\theta\, dr\, dz
\end{align*}
Thus, in particular, the vector field $A$ has only ${\bf e_\theta}$ component. Defining $\psi(r,z)$ by $A= \frac{\psi(r,z)}{r} {\bf e_\theta}$,
$\psi$ satisfies \eqref{psi_bdy} and \eqref{psi_eq} as $A$ satisfies \eqref{Aeq}.

\smallskip

If $D$ is a cylinder, the result follows by using the Green's function expansion in terms of Bessel functions. $\Box$

\medskip
Since $\psi|_{\partial D}=0$, $\nabla \psi = \psi_r\, e_r+ \psi_z \, e_z$ is normal to $\partial D$. This implies $(-\psi_z,\psi_r) \cdot n= (ru^r, ru^z)\cdot n=0$. For boundary points not on the axis, this implies $u\cdot n=0$. In the case of the ball, for the boundary points on the axis, using continuity of $u$ and the boundary, we can conclude $u^z=0$ at these points.

\medskip
In addition to showing existence of a stream function, we can use well-known results for the Poisson equation 
and Lemma $\ref{wang}$ to conclude that $\ds \frac{\psi(r,z)}{r}$ satisfies axis conditions and has regularity estimates 
in $(r,z)$ coordinates. The next theorem allows us to do this.

\smallskip
\begin{lemma}
\label{globalreg}
Let $f\in C^{k,\alpha}(\ol{D})$ and $g \in C^{k+2,\alpha}(\partial D)$ and suppose $g$, $f$ are axially symmetric. Additionally, in cylindrical coordinates $(r,z)$, suppose $f$ and $g$ satisfy
\begin{align*}
\partial_r^{2m} f(0^+,z) &=0, \quad 0\le 2m\le k \\
\partial_r^{2\ell } g (0^+,z) &=0, \quad 0\le 2\ell \le k+2
\end{align*}
Then there exists a unique $\phi(r,z)$ satisfying
\begin{align*}
\mathscr{L}\phi &= f \\
\phi|_{\partial D} &= g.
\end{align*}
In particular, $\phi$ satisfies the following estimate
\begin{align}
\label{g_est}
\|\phi\|_{C_{r,z}^{k+2,\alpha}(\ol{D})} \le C(\|g\|_{C^{k+2,\alpha}(\partial D)} + \|f\|_{C^{k,\alpha}(\ol{D})})
\end{align}
and $\phi$ satisfies
\begin{align}
\label{phipole}
\partial_r^{2m} \phi(0^+,z) &=0, \quad 0\le 2m\le k+2.
\end{align}
Suppose $D' \subset\subset D$ and $d\le \mbox{dist}(D', \partial D)$. Then we have the following interior estimate
\begin{align}
\label{interior_est}
d \|\nabla \phi\|_{L^\infty(D')} + d^2 \| \nabla^2 \phi\|_{L^\infty(D')} \le C (\|\phi\|_{L^\infty(D)} + \|f\|_{C^{\alpha}(D)})
\end{align}
\end{lemma}

\Pf Consider the following system
\begin{align*}
-\Delta \Phi &= f\, {\bf e^\theta} \\
\Phi|_{\partial D} &= g\, {\bf e^\theta}
\end{align*}
By Lemma \ref{wang}, $f\, {\bf e^\theta}$ corresponds to a $C^{k,\alpha}$ vector field in Cartesian coordinates and similarly $g\, {\bf e^\theta}$ corresponds to a $C^{k+2,\alpha}$ vector field. Then by well-known results for the Poisson equation (see \cite{Kellogg} or Theorems 6.6,  6.19, and Corollary 6.3 of \cite{GilbargTrudinger}), there exists a unique vector field $\Phi$ in $C^{k+2,\alpha}(\ol{D})$ that satisfies the above system. Additionally, $\Phi = \phi(r,z)\, {\bf e^\theta}$
for some $\phi$. Applying Lemma $\ref{wang}$, we get that $\phi\in C_{r,z}^{k+2}(\ol{D})$ and satisfies $\eqref{phipole}$. In addition by restricting $\Phi$ to the $\theta=0$ plane, the fact that $\Phi$ is H\"{o}lder continuous 
implies that $\phi$ is also H\"{o}lder continuous so $\phi \in C_{r,z}^{k+2,\alpha}(\ol{D})$. By construction, $\mathscr{L} \phi=f$. The derivative estimates, \eqref{g_est} and \eqref{interior_est}, for $\phi$ follow from the analogous classical estimates (Corollary 6.3 and Theorem 6.19 of \cite{GilbargTrudinger}) for solutions of the Poisson equation applied to $\Phi$. 

$\Box$

\bigskip
By applying the above lemma to the function $\ds \frac{\psi(r,z)}{r}$, \eqref{phipole} gives us the desired axis conditions for the stream function. With the two lemmas above at our disposal, we can continue our calculation of the stream function $\psi$.

Using \eqref{Apsi}, we can write the stream function $\psi$ as
\begin{align*}
\psi(\bar{r},\bar{z}) = \int_D \bar{r}\mathscr{A}^\theta(\bar{r},\bar{z},r,z)\omega^\theta(r, z)\, dr\, dz
\end{align*}
Let $x=\frac{\bar{r}-r}{r}$ and $y=\frac{\bar{z}-z}{r}$. We can write the integral corresponding to the first term in \eqref{eq:A} as
\begin{align*}
\frac{1}{2\pi\sqrt{1+x}}\int_0^\pi \frac{\cos \theta\, d\theta}{\sqrt{2(1-\cos \theta)+\frac{x^2+y^2}{1+x}}}
\end{align*}
Define 
\begin{align}\label{Fs}
F(s)=\int_0^\pi \frac{\cos\theta\, d\theta}{\sqrt{2(1-\cos\theta)+s}}.
\end{align}
The function $F$ cannot be expressed in terms of elementary functions. A formula for $F$ in terms of elliptic integrals can be found in Lamb \cite{Lamb}.
After some computation, $ \frac{1}{\sqrt{1+x}} = \sqrt{\frac{r}{\bar{r}}}$ and $ \frac{x^2+y^2}{1+x} = \frac{(r-\bar{r})^2+(z-\bar{z})^2}{r\bar{r}}$.
Then the stream function  $\psi$ can be written as
\begin{align*}
\psi(\bar{r},\bar{z}) &=\frac{1}{2\pi}\int_D\sqrt{r\bar{r}}F\left(\frac{(r-\bar{r})^2+(z-\bar{z})^2}{r\bar{r}}\right)\omega^{\theta}(r,z)dr\, dz\\
&\, + \frac{1}{4\pi}\int_D \int_0^{2\pi} r\bar{r} \cos\theta \cdot h(\bar{r},\bar{z},r,\theta,z)\omega^{\theta}(r,z)\, d\theta dr\, dz.
\end{align*}
Define 
\begin{align*}
H(\bar{r},\bar{z},r,z)=\frac{1}{4\pi}\int_0^{2\pi} r\bar{r} \cos\theta \cdot h(\bar{r},\bar{z},r,\theta,z) d\theta.
\end{align*}
We can compute and integrate by parts to see that
\begin{align*}
-H_{zz}+\frac{H_r}{r} -H_{rr} &= \frac{\bar{r}r}{4\pi} \int_0^{2\pi} \left(-h_{rr}-\frac{h_r}{r}-h_{zz}+\frac{h}{r^2}\right)\cos\theta\, d\theta \\
&= \frac{\bar{r}r}{4\pi}  \int_0^{2\pi} \left(-h_{rr}-\frac{h_r}{r}-h_{zz}-\frac{h_{\theta\theta}}{r^2}\right)\cos\theta\, d\theta \\
&= -\frac{\bar{r}r}{4\pi}  \int_0^{2\pi} (\Delta h)(\bar{r},\bar{z},r,\theta,z)\, \cos\theta\, d\theta=0.
\end{align*}

\medskip
We summarize the above calculations with the following proposition.

\begin{proposition}
\label{GREEN}
The function 
\begin{align}
\mathscr{G}(\bar{r},\bar{z},r,z) = \frac{\sqrt{r\bar{r}}}{2\pi}F\left(\frac{(r-\bar{r})^2+(z-\bar{z})^2}{r\bar{r}}\right)+H(\bar{r},\bar{z},r,z).
\end{align}
is the Green's function for the operator $L$ for the domain $D$. The function $H$ satisfies the following:
\begin{align}
\mathscr{L} \left(\frac{H}{r}\right)&= r\cdot LH=-\frac{H_{zz}}{r} + \frac{H_r}{r^2}-\frac{H_{rr}}{r} = 0 \label{Hcorrector}, \quad r>0\\
H(\bar{r},\bar{z},r,z) &= - \frac{\sqrt{r\bar{r}}}{2\pi} 
F\left( \frac{(r-\bar{r})^2+(z-\bar{z})^2}{r\bar{r}}\right)
, \quad (r,z)\in \partial D \nonumber 
\end{align}
with the axis condition
\begin{align*}
\partial_r^{(2m)}\left.\left( \frac{H(\bar{r},\bar{z},r,z)}{r}\right)\right|_{r=0^+}=0, \quad m=0,1,2,\ldots
\end{align*}

In addition, we have 

\begin{align*}
u(\bar{r},\bar{z})= \frac{1}{\bar{r}}\int_D \nabla_{(\bar{r},\bar{z})}^\perp \mathscr{G}(\bar{r},\bar{z},r,z) \omega^\theta(r,z) \, dr\, dz = \frac{1}{\bar{r}}\int_D r\nabla_{(\bar{r},\bar{z})}^\perp \mathscr{G}(\bar{r},\bar{z},r,z) w(r,z) \, dr\, dz
\end{align*}
\end{proposition}

\bigskip
\noindent
\begin{remark}
As seen by the above calculations, the Green's function for $L$ is related to the Green's function $G_D$ for the 3D Laplacian by
\begin{align*}
\mathscr{G}(\bar{r},\bar{z},r,z) = \int_0^{2\pi} r\bar{r} \cdot G_D(\bar{r},\bar{z},r,z)\cos \theta d\theta.
\end{align*}
\end{remark}

\bigskip
After some computation, we have the following expressions for $u^r$ and $u^z$ on $D$:
\begin{align}
\label{ur}
u^r(\bar{r},\bar{z}) &=  \int_D \left[\frac{(z-\bar{z})\sqrt{r}}{\pi \bar{r}^{3/2}} F'\left(\frac{(r-\bar{r})^2+(z-\bar{z})^2}{\bar{r}r}\right) +\frac{r}{\bar{r}}\partial_{\bar{z}} H\right] w(r,z)\, dr\, dz\\
u^z(\bar{r},\bar{z}) &= \int_D\left[ \wtilde{J}(\bar{r},\bar{z}, r,z) +\frac{r}{\bar{r}}\partial_{\bar{r}} H\right] w(r,z)\, dr\, dz
\end{align}
where we define
\begin{align}
\label{scriptJ}
\wtilde{J}(\bar{r},\bar{z},r,z) = \left(\frac{r}{\bar{r}}\right)^{3/2} \wtilde{\mathscr{J}}(\bar{r},\bar{z},r,z)
\end{align}
and the function $\wtilde{\mathscr{J}}(\bar{r},\bar{z},r,z)$ is given by
\begin{align*}
\wtilde{\mathscr{J}}(\bar{r},\bar{z},r,z)= \frac{1}{\pi}\frac{(\bar{r}-r)}{r}F'\left(\frac{(r-\bar{r})^2+(z-\bar{z})^2}{\bar{r}r}\right) &\,&\\
 +\frac{1}{4\pi}\left[ F\left(\frac{(r-\bar{r})^2+(z-\bar{z})^2}{\bar{r}r}\right)-2\frac{(r-\bar{r})^2+(z-\bar{z})^2}{\bar{r}r}F'\left(\frac{(r-\bar{r})^2+(z-\bar{z})^2}{\bar{r}r}\right)\right]. 
\end{align*}
In addition, we set
\begin{align}
\label{scriptK}
\wtilde{K}(r,\bar{r},z,\bar{z}) = \frac{(z-\bar{z})\sqrt{r}}{\pi \bar{r}^{3/2}} F'\left(\frac{(r-\bar{r})^2+(z-\bar{z})^2}{\bar{r}r}\right).
\end{align}

\medskip
For later use, we also define the kernels in the integrals for $u^r$ and $u^z$
with the corrector
\begin{align}
\label{K_def}
K(\bar{r},\bar{z},r,z) := -\frac{r}{\bar{r}}\partial_{\bar{z}}\mathscr{G}  = \wtilde{K}(r,\bar{r},z,\bar{z}) +\frac{r}{\bar{r}}\partial_{\bar{z}} H
\end{align}
\begin{align}
\label{L_def}
J(\bar{r},\bar{z},r,z) := \frac{r}{\bar{r}} \partial_{\bar{r}}\mathscr{G} = \wtilde{J}(\bar{r},\bar{z},r,z) + \frac{r}{\bar{r}}\partial_{\bar{r}} H.
\end{align}

\subsection{Behavior of $F$}

\medskip
We derive estimates for the function $F$ defined in \eqref{Fs} that will be used frequently later. The details for these estimates can be found in the appendix. 

\medskip
As $F$ is tied with the Green's function of $L$ from $\eqref{GREEN}$, one may expect that $F$ behaves
 roughly logarithmically. However, $F$ will have better asymptotic properties than $\log$ in certain regimes. This will be
key for our estimates later.

First, one can bound $F$ easily with
\begin{align}
|F(s)| \lesssim \left(\frac{1}{s}\right)^{1/2}
\end{align}
but in fact, we have even better asymptotics at $s=0$ and $s=\infty$
\begin{align}
\label{Fexp}
F(s) &= -\frac{1}{2}\log(s)+\log 8-2+O(s\, \log(s))\quad \mbox{near} \quad s=0 \\
F(s) &= \frac{\pi}{2} \frac{1}{s^{3/2}} + O(s^{-5/2}) \quad \mbox{near} \quad s=\infty 
\end{align}
and expansions gotten by  formally differentiating the series holds. They are as follows:
\begin{align} 
\label{F'exp0}
F'(s) &=-\frac{1}{2}\frac{1}{s}+O(\log s) \quad \mbox{near} \quad s=0\\
\label{F''exp0}
F''(s) &=\frac{1}{2} \frac{1}{s^2}+O(1/s)\quad \mbox{near} \quad s=0
\end{align}
\begin{align}
\label{F'expinfty}
F'(s) &= -\frac{3\pi}{4} \frac{1}{s^{5/2}} + O(s^{-7/2}) \quad \mbox{near} \quad s=\infty \\
F''(s) &= \frac{15\pi}{8} \frac{1}{s^{7/2}} + O(s^{-9/2}) \quad \mbox{near} \quad s=\infty.
\end{align}
Let $\epsilon>0$ be a small constant such that for $0<s<\epsilon$, all the expansions above for $F,F'$, and $F''$ near $0$ are valid. We will refer to this $\epsilon$ in our proofs later.

\medskip
We summarize upper bounds on $F$ in the following lemma which is from \cite{FengSverak}
\begin{lemma}
\label{Fest}
For every non-negative integer $k$, for all $s>0$
\begin{align}
|F(s)| & \lesssim_\tau \min \left( \left(\frac{1}{s}\right)^\tau, \left(\frac{1}{s}\right)^{3/2}\right), \quad 0<\tau\le\frac{1}{2} \label{Fbound} \\
|F^{(k)}(s)| &\lesssim_k  \min \left( \left(\frac{1}{s}\right)^k, \left(\frac{1}{s}\right)^{k+3/2}\right), \quad k>0. \label{Fprimebound}
\end{align}
\end{lemma}
We will use the above bounds constantly throughout the rest of our proofs.

\section{Gradient upper bound}

Our first goal is to prove a Kato type estimate on $\|\nabla u\|_\infty$ (see \cite{Kato}), which will imply an upper bound of $\|\nabla w\|_\infty$. Our estimate will have parallels with the analogous estimate for $\|\nabla u\|_\infty$ for the  2D Euler equations, but the estimates become more tedious due to the more complex Biot-Savart law.

\subsection{Some Green's function computations and derivative estimates}

\medskip
\noindent
 Here, we will collect computations concerning the kernels $\wtilde{K}$ and $\wtilde{J}$ arising in the integrals for $u^r$ and $u^z$ which will be useful in our later estimates. 
In order to have estimates on $\nabla u$ we will need to bound derivatives of $\wtilde{K}$ and $\wtilde{J}$ which are 
\begin{align}
\partial_{\bar{r}}\wtilde{K} (r,\bar{r},z,\bar{z}) &=-\frac{3}{2\pi} \frac{(z-\bar{z}) \sqrt{r}}{\bar{r}^{5/2}} F'\left(s\right) + \frac{z-\bar{z}}{\pi\bar{r}^{3/2}} \sqrt{r} F''\left(s\right)\left( \frac{-2\bar{r}(r-\bar{r})-((r-\bar{r})^2+(z-\bar{z})^2)}{(\bar{r})^2 r}\right) \nonumber \\
&= -\frac{(z-\bar{z})\sqrt{r}}{\pi\bar{r}^{5/2}}\left[\frac{3}{2} F'(s)+sF''(s)\right]-2 \frac{(z-\bar{z})(r-\bar{r})}{\pi\bar{r}^{5/2}\sqrt{r}}F''(s)
\label{Kr} \\
\partial_{\bar{z}} \wtilde{K}(r,\bar{r},z,\bar{z}) &= -\frac{\sqrt{r}}{\pi\bar{r}^{3/2}} \left[ F'\left(s\right)+2\frac{(z-\bar{z})^2}{\bar{r}r} F''\left(s\right) \right] \label{Kz}
\end{align}
\begin{align}
\partial_{\bar{r}} \wtilde{J}(r,\bar{r},z,\bar{z}) &= \frac{1}{\pi} \left(\frac{\sqrt{r}}{\bar{r}^{3/2}}- \frac{3}{2}\frac{(\bar{r}-r)\sqrt{r}}{\bar{r}^{5/2}}\right)F'(s)+ \frac{1}{\pi} \frac{(\bar{r}-r)\sqrt{r}}{\bar{r}^{3/2}} F''(s) (\partial_{\bar{r}}s) \nonumber  \\
&\,\,\, -\frac{3}{8\pi}\frac{r^{3/2}}{\bar{r}^{5/2}} \left[ F\left(s \right)-2s F'\left(s \right)\right]+ \frac{1}{4\pi}\left[-F'(s) (\partial_{\bar{r}} s)-2s F''(s) (\partial_{\bar{r}} s)\right] \left(\frac{r}{\bar{r}}\right)^{3/2} \nonumber \\
&= \frac{1}{\pi} \left(\frac{\sqrt{r}}{\bar{r}^{3/2}}- \frac{3}{2}\frac{(\bar{r}-r)\sqrt{r}}{\bar{r}^{5/2}}\right)F'(s)-
\frac{1}{\pi} \frac{(\bar{r}-r)\sqrt{r}}{\bar{r}^{3/2}} F''(s) \left(\frac{2(r-\bar{r})}{\bar{r}r} + \frac{s}{\bar{r}}\right)
 \nonumber \\
 &\,\,\, -\frac{3}{8\pi}\frac{r^{3/2}}{\bar{r}^{5/2}} \left[ F\left(s \right)-2s F'\left(s \right)\right]+ \frac{1}{4\pi}
 [F'(s)+2sF''(s)]\left(\frac{r}{\bar{r}}\right)^{3/2}\left(\frac{2(r-\bar{r})}{\bar{r}r} + \frac{s}{\bar{r}}\right) \\
&= \frac{1}{\pi} \frac{r^{3/2}}{\bar{r}^{5/2}}\left[- \frac{3}{8}F(s)
+ sF'(s)+\frac{1}{2}s^2 F''(s)\right]+
\frac{2(r-\bar{r})}{\pi} \frac{\sqrt{r}}{\bar{r}^{5/2}}[F'(s)+sF''(s)] \label{Jr}\\
&\,\,+\frac{1}{\pi}\frac{\sqrt{r}}{\bar{r}^{3/2}}F'(s)+
\frac{2}{\pi}\frac{(r-\bar{r})^2}{\bar{r}^{5/2}\sqrt{r}}F''(s) \\
\partial_{\bar{z}} \wtilde{J}(r,\bar{r},z,\bar{z}) &= \frac{1}{\pi}\frac{(\bar{r}-r)\sqrt{r}}{\bar{r}^{3/2}} F''(s) \left(-\frac{2(z-\bar{z})}{\bar{r}r}\right) + \frac{1}{4\pi} \left[-F'(s)-2s F''(s)\right]\left(\frac{r^{3/2}}{\bar{r}^{3/2}}\right)\left(-\frac{2(z-\bar{z})}{\bar{r}r}\right) \label{Jz}
\end{align}
where we have defined $\ds s =\frac{(r-\bar{r})^2+(z-\bar{z})^2}{\bar{r}r}$ and $\ds \partial_{\bar{r}} s= -\frac{2(r-\bar{r})}{\bar{r}r} - \frac{s}{\bar{r}}$.

\medskip
\bigskip
\noindent
For later use, we define 
\begin{align*}
x=(\bar{r},\bar{z}), \quad \mbox{and}\,\,\, y=(r,z).
\end{align*}

\begin{remark}
	Observe that up to factors of $r$ and $\bar{r}$, the most singular terms for the derivatives of $\wtilde{K}$
	and $\wtilde{J}$ above  are similar to those gotten by computing the second derivatives of the 2D Laplacian Green's function. This will lead
	to the double exponential upper bound in Theorem \ref{upper}
\end{remark}

\medskip
Additionally, we will need bounds for the Green's function $\mathscr{G}$ of $L$, which will then allow for bounds on derivatives of $J$ and $K$.

\begin{proposition}
\label{Hest}
The function $\mathscr{G}$ defined by Proposition \ref{GREEN} satisfies the following estimates for $\bar{r}>0$

\begin{align*}
\left| \nabla_{\bar{r},\bar{z}}^2 \left(\frac{\mathscr{G}(\bar{r},\bar{z},r,z)}{\bar{r}}\right)\right| & \le C(D)\min \left( \frac{r}{|x-y|^3}, \sqrt{\frac{r}{\bar{r}}}\frac{1}{|x-y|^2}\right).
\end{align*}
\end{proposition}

\Pf

Recall
\begin{align*}
\frac{\mathscr{G}(\bar{r},\bar{z},r,z)}{\bar{r}} = \int_0^{2\pi} r\cos \theta \cdot G_D(\bar{r},\bar{z},r,\theta,z)\, d\theta
\end{align*}
where $G_D$ is defined through $\eqref{3DGreen}$. Using the classical Green's function bound \cite{Widman, Krasovskii}
\begin{align*}
|\nabla^2 G_D(x,y)| \le \frac{C(D)}{|x-y|^3},
\end{align*}
we can arrive at the bound
\begin{align*}
 \left|\nabla_{\bar{r},\bar{z}}^2 \left(\frac{\mathscr{G}}{\bar{r}}\right)\right| & \le C(D) \int_0^{2\pi} \frac{r}{(r^2-2r\bar{r}\cos\theta+\bar{r}^2+(z-\bar{z})^2)^{3/2}}\, d\theta \\
&= \frac{C(D)}{r^{1/2}\bar{r}^{3/2}} \int_0^{2\pi} \frac{1}{\left(2(1-\cos \theta) + \frac{(r-\bar{r})^2+(z-\bar{z})^2}{r\bar{r}}\right)^{3/2}}\, d\theta.
\end{align*}
As before, we set $\ds s= \frac{(r-\bar{r})^2+(z-\bar{z})^2}{r\bar{r}}$. Then applying a similar argument as in section {\bf A.1} of the Appendix we can arrive at 
\begin{align}
\label{s_bound}
 \int_0^{2\pi} \frac{d\theta}{\left(2(1-\cos \theta) +s\right)^{3/2}} \lesssim \min\left(\frac{1}{s}, \frac{1}{s^{3/2}}\right).
\end{align}
Indeed, easily we have 
\begin{align*}
 \int_0^{2\pi} \frac{d\theta}{\left(2(1-\cos \theta) +s\right)^{3/2}} \le \frac{2\pi}{s^{3/2}}.
\end{align*}
For the other possible upper bound can rewrite the integral as
\begin{align*}
 \int_0^{2\pi} \frac{d\theta}{\left(2(1-\cos \theta) +s\right)^{3/2}} = \frac{1}{2} \int_0^{\pi/2} \frac{1}{(\sin^2 \varphi+s/4)^{3/2}}\, d\varphi
\end{align*}
 and then by \eqref{Fprime1} we get an upper bound of a constant times $1/s$.
Using the bound \eqref{s_bound}, we arrive at the desired estimate for $\nabla^2( \mathscr{G}/\bar{r})$. $\Box$

\subsection{Kato estimate}

Next, we prove the key estimate that will allow us to deduce Theorem \ref{upper}. For the 2D Euler equations, this type of estimate was proven by Kato \cite{Kato}
\begin{theorem} (Kato type estimate)
\label{kato}
Let $w\in C^\alpha(D)$, $\alpha>0$. Fix $R>0$ and let $D_R=\{(r,z): (r,z)\in D\,\,\mbox{and}\,\, r <R\}$.
\begin{align}
\label{kato_est}
\|\nabla u\|_{L^\infty(D_R)} \le C_1(\alpha,D)\|w_0\|_\infty\left(1+(R+R^3)\log\left(1+\frac{\|w\|_{C^\alpha}}{\|w_0\|_\infty}\right)\right)
\end{align}

\end{theorem}

\begin{remark}
Observe that as we are closer to the axis, the effect of the logarithm on the right hand side is diminished. It is in this respect that the above estimate is different than the estimate for $\|\nabla u\|_\infty$ for solutions of 2D Euler.
\end{remark}

\subsection{Proof of Theorem \ref{kato}}

Let $x=(\bar{r},\bar{z})$. Define  $\epsilon_0= ((\sqrt{5}-1)/4)\epsilon $. Recall $\epsilon$ is the radius of the ball centered at $0$ for which our expansions for $F$ in section 2.3 hold. Observe that
\begin{align*}
B_{\epsilon_0\bar{r}}(x) \subset \{ (r,z): s(r,z) <\epsilon\}
\end{align*}
so we can use the expansion for $F(s)$ close to $s=0$ on this ball later.
Let 
\begin{align}
\label{delta_choice}
\delta =\min \left(c, \epsilon_0/2, \left( \frac{\|w_0\|_{L^\infty}}{\|w\|_{C^\alpha}}\right)^{1/\alpha}\right)
\end{align}
where $c$ is chosen small so that the set of $x$ such that $\mbox{dist}(x,\partial D)> 2\delta$ is non-empty.
We will bound
\begin{align*}
\partial_{\bar{r}}u^r(\bar{r},\bar{z})= \int_D\partial_{\bar{r}}K(r,z,\bar{r},\bar{z})w(r,z)\, dr\, dz.
\end{align*}
Proposition \ref{Hest} will allow for better decay estimates for the kernel when $|x-y|$ is large
which will lead to the decay factor in front of the logarithm term in $\eqref{kato_est}$.

By the incompressibility condition and axis condition, bounding $\partial_{\bar{r}} u^r$ will imply the desired
 bound on $\partial_{\bar{z}} u^z$. We will sketch the proof for the derivatives $\partial_{\bar{r}} u^z$ and $\partial_{\bar{z}} u^r$ in the appendix.

\bigskip
\noindent
{\bf Case 1: $\mbox{dist}(x,\partial D) > 2\delta$}

\medskip
We divide the integral for $\partial_{\bar{r}}u^r$ into three regions
\begin{align*}
\partial_{\bar{r}}u^r(\bar{r},\bar{z}) &=\left( \int_{B_{\delta \bar{r}}(x)} +\int_{\Omega \cap B_{\delta\bar{r}}^c(x)} +\int_{D\setminus \Omega}\right)\partial_{\bar{r}}K(r,z,\bar{r},\bar{z})w(r,z)\, dr\, dz \\ 
&= I+II+III
\end{align*}
where $\Omega = \{(r,z)\in D:  \frac{1}{2}\bar{r} < r < 2\bar{r}, 0\le z\le 1\}$. 
Recall from above that we have defined
\begin{align*}
\partial_{\bar{r}}K(r,z,\bar{r},\bar{z}) = \partial_{\bar{r}}\wtilde{K}(r,z,\bar{r},\bar{z}) + \partial_{\bar{r}}\left( \frac{r}{\bar{r}} \partial_{\bar{z}} H\right).
\end{align*}
For $I$, we can use the expansion \eqref{Kr} for $\partial_r \wtilde{K}$ and \eqref{F''exp0} to get that the most dangerous term is
\begin{align}
\label{danger}
-\frac{1}{\pi\sqrt{\bar{r}}} \int_{B_{\delta\bar{r}}(x)}  \frac{r^{3/2}(z-\bar{z})(r-\bar{r})}{|x-y|^4} w(r,z)\, dr\, dz.
\end{align}
Indeed, using $\eqref{Kr}$, $\eqref{F'exp0}$, and $\eqref{Fprimebound}$ we can bound the remainder terms
\begin{align*}
\left|\partial_{\bar{r}} \wtilde{K} + \frac{1}{\pi\sqrt{\bar{r}}}\frac{r^{3/2}(z-\bar{z})(r-\bar{r})}{|x-y|^4}\right| \lesssim \left( \frac{|z-\bar{z}|\sqrt{r}}{\bar{r}^{5/2}}\right)\left(\frac{1}{s}\right) \lesssim \frac{r^{3/2}}{\bar{r}^{3/2}} \frac{1}{|x-y|}.
\end{align*}
Thus, the terms other than $\eqref{danger}$ from the expansion of $\partial_r \wtilde{K}$ can be controlled since
\begin{align*}
\frac{1}{\bar{r}^{3/2}}\int_{B_{\delta\bar{r}}(x)} \frac{r^{3/2}}{|x-y|}w(y,t)\, dy \le C\bar{r}^{2} \|w_0\|_\infty.
\end{align*}

By writing the kernel in the integral $\eqref{danger}$ in polar coordinates $\rho, \phi$ centered at $x=(\bar{r},\bar{z})$, we get
\begin{align}
\label{taylor}
- \frac{(\bar{r}+\rho \cos\phi)^{3/2} \cos \phi\sin\phi}{\rho^2} = - \frac{\bar{r}^{3/2} \cos\phi\sin\phi}{\rho^2} + \bar{r}^{1/2} O(1/\sqrt{\rho}).
\end{align}
When integrated over $B_{\delta\bar{r}}(x)$, the second term in \eqref{taylor} is controlled by $C\bar{r}^{2} \|w_0\|_{L^\infty}$. For the other term,
\begin{align*}
\left| \frac{1}{\sqrt{\bar{r}}} \int_{B_{\delta\bar{r}}(x)}  \frac{\bar{r}^{3/2}\cos\phi\sin\phi}{\rho^2} w(\rho,\phi) \rho\, d\rho\, d\phi\right|  &= \left|\bar{r} \int_{B_{\delta\bar{r}}(x)} \frac{\cos\phi\sin\phi}{\rho} (w(\rho, \phi)-w(\bar{r},\bar{z}))\, d\rho\, d\phi\right| \\
& \le C\bar{r}\|w\|_{C^\alpha} \int_0^{\delta\bar{r}} \rho^{-1+\alpha}\, d\rho \\ &\le C\bar{r}^{1+\alpha} \delta^\alpha \|w\|_{C^\alpha} \le C(\alpha)\bar{r}^{1+\alpha} \|w_0\|_{L^\infty}.
\end{align*}
Now, we estimate the contribution from the corrector function $H$.  For $y=(r,z)\in B_{\delta\bar{r}}(x)$, the function $H$ satisfies

\begin{align*}
\mathscr{L}_{\bar{r},\bar{z}} \left( \frac{H(\bar{r}, \bar{z}, r, z)}{\bar{r}}\right) &= 0 \quad\quad\mbox{for} \quad (\bar{r},\bar{z})\in D \\
\frac{H(\bar{r}, \bar{z}, r, z)}{\bar{r}} &= - \frac{1}{2\pi} \sqrt{\frac{r}{\bar{r}}} F\left( \frac{(r-\bar{r})^2+(z-\bar{z})^2}{r\bar{r}}\right) \quad\quad\mbox{for} \quad (\bar{r},\bar{z}) \in \partial D.
\end{align*}
By using bounds for $F$ we obtain,
\begin{align*}
\sup_{(\wtilde{r},\wtilde{z})\in \partial D}  \sqrt{\frac{r}{\wtilde{r}}}\left| F\left( \frac{(r-\wtilde{r})^2+(z-\wtilde{z})^2}{r\wtilde{r}}\right)\right| &\le
 \sup_{(\wtilde{r},\wtilde{z})\in \partial D, 0\le \wtilde{r}\le \delta} \left|\sqrt{\frac{r}{\wtilde{r}}} F\left( \frac{(r-\wtilde{r})^2+(z-\wtilde{z})^2}{r\wtilde{r}}\right)\right| \\ 
&\quad +  
\sup_{(\wtilde{r},\wtilde{z})\in \partial D,  \wtilde{r}\ge \delta} \left|\sqrt{\frac{r}{\wtilde{r}}} F\left( \frac{(r-\wtilde{r})^2+(z-\wtilde{z})^2}{r\wtilde{r}}\right)\right| \\
& \le  C\sup_{(\wtilde{r},\wtilde{z})\in \partial D, 0\le \wtilde{r}\le \delta} \left|\frac{\wtilde{r}r^2}{((r-\wtilde{r})^2+(z-\wtilde{z})^2)^{3/2}}\right| \\
& \quad + C \sup_{(\wtilde{r},\wtilde{z})\in \partial D,  \wtilde{r}\ge \delta}  \sqrt{\frac{r}{\wtilde{r}}} \left| \log\left( \frac{\wtilde{r}r}{(r-\wtilde{r})^2+(z-\wtilde{z})^2}\right) \right|\\
& \le C+C\log \delta^{-1} \le C\log \delta^{-1}
\end{align*}
where we can make $\delta$ smaller if needed.

By relating this equation to the Poisson equation as in the proof of Lemma 2.5, we can use the maximum principle so that for $\ds y\in B_{\delta\bar{r}}(x)$
\begin{align*}
\left| \frac{H(\bar{r},\bar{z},r,z)}{\bar{r}}\right| \le \sup_{(\wtilde{r},\wtilde{z})\in \partial D}  \sqrt{\frac{r}{\wtilde{r}}}\left| F\left( \frac{(r-\wtilde{r})^2+(z-\wtilde{z})^2}{r\wtilde{r}}\right)\right| \le C \log \delta^{-1}.
\end{align*}
Using the interior estimate \eqref{interior_est} of Lemma \ref{globalreg}, we get
\begin{align*}
\left|\nabla_{\bar{r},\bar{z}}^2 \left(\frac{H(\bar{r},\bar{z},r,z)}{\bar{r}}\right)\right|
\le C(D)  \frac{1}{\delta^2}\sup_{(\bar{r},\bar{z})\in D}\left|\frac{H(\bar{r},\bar{z},r,z )}{\bar{r}}\right|
 \le C(D)  \frac{1}{\delta^2}\log \delta^{-1}.
\end{align*}
Inserting this estimate to the integral of the corrector over $B_{\delta\bar{r}}(x)$:
\begin{align*}
\left|\int_{B_{\delta\bar{r}}(x)} r \partial_{\bar{r}}
\partial_{\bar{z}} \left( \frac{H}{\bar{r}}\right) w(r,z)\, dr\, dz\right|   \le C\bar{r}^3 \|w_0\|_{L^\infty} \log \delta^{-1}.
\end{align*}
Combining estimates for $I$, we get 
\begin{align}
\label{I_integral}
|I| \le C(\alpha) (\bar{r}^{2}+\bar{r}^3\log \delta^{-1}+\bar{r}^{1+\alpha}) \|w_0\|_\infty. 
\end{align}

Now, we will bound $II$ and we use
\begin{align*}
\partial_{\bar{r}}K(r,z,\bar{r},\bar{z})=-r \partial_{\bar{r}}\partial_{\bar{z}} \left(\frac{\mathscr{G}}{\bar{r}}\right) .
\end{align*}

\medskip
\noindent
Using Proposition \ref{Hest},
\begin{align}
|II|\le \int_{\Omega \cap B_{\delta\bar{r}}^c(x)} \left| \partial_{\bar{r}}K(r,z,\bar{r},\bar{z}) w(r,z)\right|\, dr\, dz &\le C(D)\|w_0\|_\infty  \int_{\Omega \cap B_{\delta\bar{r}}^c(x)} \frac{r^{3/2}}{\bar{r}^{1/2}} \frac{1}{|x-y|^2}\, dr\, dz \nonumber \\
&\le C(D)\bar{r} \|w_0\|_\infty \int_{\Omega \cap B_{\delta\bar{r}}^c(x)}\frac{1}{|x-y|^2}\, dr\, dz 
\nonumber \\  & \le C(D) \bar{r} \|w_0\|_\infty (1+ \log \bar{r}^{-1} +\log \delta^{-1}) \label{II_integral}
\end{align}

\bigskip

\noindent
Using the other bound from Proposition \ref{Hest}, we can bound $III$ 

\begin{align}
|III| &\le \int_{D\setminus\Omega } \left|\partial_{\bar{r}}K(r,z,\bar{r},\bar{z}) w(r,z)\right|\, dr\, dz\le C(D) \|w_0\|_\infty  \int_{D\setminus\Omega } \frac{r^2}{|x-y|^3}\, dr\, dz \nonumber \\
&\le C(D) \|w_0\|_\infty   \int_{D\setminus\Omega } \frac{1}{|x-y|}\, dr\, dz  \le C(D) \|w_0\|_\infty. \label{III_integral}
\end{align}
Above, we have used that $r^2 \le C(r-\bar{r})^2 $ for some constant $C$ on our domain of integration. 
After combining the estimates \eqref{I_integral}, \eqref{II_integral}, and \eqref{III_integral}, we get the desired estimates at interior points.

\medskip
{\bf Case 2: $\mbox{dist} (x,\partial D) < 2\delta$}

\medskip
We can express the derivatives of $u^r$ and $u^z$ as
\begin{align*}
\partial_{\bar{r}} u^r &= -\partial_{\bar{r}}\partial_{\bar{z}} \left( \frac{\psi}{\bar{r}}\right) \\
\partial_{\bar{z}} u^r &= - \partial_{\bar{z}}^2  \left( \frac{\psi}{\bar{r}}\right) \\
\partial_{\bar{r}} u^z &= \partial_{\bar{r}} \left(\frac{\partial_{\bar{r}} \psi}{\bar{r}}\right) = -\partial_{\bar{z}}^2 \left(\frac{\psi}{\bar{r}}\right)- \omega^\theta
\end{align*}
where in the last equality we use $\mathscr{L}(\psi/\bar{r})=\omega^\theta$.

\medskip
Find a point $x'$ such  that $\mbox{dist}(x,\partial D) \ge 2\delta$ and $|x'-x|\le  C(D) \delta$. By estimate \eqref{g_est} of Lemma \ref{globalreg} applied to $\frac{\psi}{r}$, we know that
\begin{align*}
|\nabla u^r(x)-\nabla u^r(x')| & \le C(\alpha ,D)\delta^\alpha \|\omega^\theta\|_{C^\alpha} \le C(\alpha ,D)\delta^\alpha \|w\|_{C^\alpha} \\
|\partial_{\bar{r}} u^z(x)-\partial_{\bar{r}} u^z(x')|& \le C(\alpha ,D)\delta^\alpha \|w\|_{C^\alpha} .
\end{align*}
Combining these above estimates with our interior estimates above\eqref{I_integral}, \eqref{II_integral}, and \eqref{III_integral} , we can deduce the main estimate $\eqref{kato_est}$. Recall our choice of $\delta$ in $\eqref{delta_choice}$. Then the two estimates we just derived
above become part of the first term on the right side of $\eqref{kato_est}$. The log factor in $\eqref{kato_est}$ will arise
from our earlier estimates $\eqref{I_integral}$, $\eqref{II_integral}$, and choice of $\delta$. We can deduce the $R$ factors in front of the log term of \eqref{kato_est} since our interior
 estimate decays with $R$ and our boundary estimate has no $\log $ term. This completes the proof of Theorem \ref{kato} $\Box$.

\subsection{Proof of Theorem \ref{upper}}

With Theorem \ref{kato} now at our disposal, we can derive our first main result. Using \eqref{kato_est}, the proof of estimate \eqref{upper1} is standard and we refer readers to \cite{KiselevSverak} for the details. To prove part (b) of the theorem, using that $u^r(0^+,z)=0$ for all $z$ and \eqref{kato_est},
\begin{align*}
\frac{d}{dt} \Phi_t^r(x) & = u^r(\Phi_t^r(x),t) \ge -\|\nabla u\|_{L^\infty(D_{2\Phi_t(x)})}\Phi_t^r(x) \\
& \ge -C\left(1+ \Phi_t^r(x) \log \left(1+ \frac{\|\nabla w\|_{\infty}}{\|w_0\|_\infty}\right)\right) \Phi_t^r(x) \\
& \ge -C\left(1+ \Phi_t^r(x) \exp(C\|w_0\|_\infty t)\right) \Phi_t^r(x).
\end{align*}
In the last inequality above, we used the double exponential upper bound from part (a). From the above differential inequality,
it is not hard to deduce the desired estimate on $\Phi_t^r(x)$. Part (c) of the theorem follows from part (b) 
and $w\circ \Phi = w_0$. $\Box$

\section{Vorticity Gradient Growth on the Boundary Away from the Axis}

In this section, we will provide an example of double exponential gradient growth of vorticity at the boundary using the ideas of Kiselev and Sverak \cite{KiselevSverak}. Observe that away from the axis $\nabla w \approx \nabla \omega^\theta$. We will take our domain to be the unit sphere. We choose the sphere as we have an explicit expression for the Green's function of $L$ for this domain which is
\begin{align}
\label{sphere}
G(r,z,\bar{r},\bar{z})= \frac{\sqrt{r\bar{r}}}{2\pi} \left( F\left(
\frac{(r-\bar{r})^2+(z-\bar{z})^2}{r\bar{r}}\right) -
F\left( \frac{(r^\ast-\bar{r})^2+(z^\ast-\bar{z})^2}{r^\ast\bar{r}}\right)\right)
\end{align}
where we define $\ds r^\ast = \frac{r}{r^2+z^2}$ and $\ds z^\ast = \frac{z}{r^2+z^2}$. Using the methods of \cite{Xiao}, we expect Theorem \ref{lower} to hold for more general domains such as a cylinder.

\medskip
The desired growth of vorticity is achieved by establishing
a ``hyperbolic flow" scenario at the boundary of the sphere. Due to axial symmetry, it is sufficient to consider everything that follows to be on a semi-circular slice of the sphere
$$
D :=\{ (r,\theta,z): \theta=0, \quad r\ge 0, \quad r^2+z^2\le 1\}
$$
and without loss of generality, we will omit the $\theta$ component in our coordinate expressions. We will consider initial data $w_0$ which is odd with respect to $z$ and positive for $z>0$. Our goal is to show that the boundary point $e_1=(r=1,z=0)$ will act as a hyperbolic fixed point of the flow near the boundary. Because of the symmetry assumptions, we can write $u$ as
\begin{align}
u(x) &= \frac{1}{2\pi\bar{r}}\int_D r\nabla^\perp G(y,x)\omega(y,t)\, dy \nonumber \\
\label{Uball}
&= \frac{1}{2\pi\bar{r}}\int_{D^+} r\nabla^\perp \left( \sqrt{r\bar{r}} \left[F\left(\frac{|x-y|^2}{r\bar{r}}\right) -F\left( \frac{|x-y^\ast|^2}{r^\ast\bar{r}}\right)-F\left(\frac{|\wtilde{x}-y|^2}{r\bar{r}}\right)+F\left(\frac{|\wtilde{x}-y^\ast|^2}{r^\ast\bar{r}}\right)\right]\right)\omega(y)\, dy
\end{align}
where we have defined 
$$
x=(\bar{r},\bar{z}),\, y=(r,z),\, \wtilde{x}=(r,-z),\, y^\ast=(r^\ast,z^\ast),\, \mbox{and}\, D^+=D\cap \{z\ge 0\}.
$$ 

\medskip
After some computation, we have the following identities 
\begin{align}
\label{identities}
|y^\ast-e_1|^2= \frac{|y-e_1|^2}{|y|^2}, \quad \frac{z^\ast}{|y^\ast-e_1|^2} = \frac{z}{|y-e_1|^2}, \quad \frac{r^\ast-1}{|y^\ast-e_1|^2}=-1-\frac{r-1}{|y-e_1|^2}.
\end{align}

\medskip
The key observation in achieving double exponential growth is an expansion for the Biot Savart law near the fixed point. For the 2D Euler equations, this is the content of Lemma 3.1 of \cite{KiselevSverak}. We will aim to prove a similar expansion in our axially symmetric setting in the lemma below. Choose a constant $N$ such that $N<\min\{1/2,\frac{\epsilon}{8}\}$ let 
\begin{align*}
S_N &=\{ 1-N<r<1, 0<z<N\}\cap D \\
Q(\bar{r},\bar{z}) &= \{ 1-N<r<\bar{r}, \bar{z}<z<N\}\cap D.
\end{align*}
Recall we defined $\epsilon$ to be a small constant so that the expansions of $F(s)$ and its derivatives ($\eqref{Fexp}, \eqref{F'exp0}$, and $\eqref{F''exp0}$) can be used for $0<s<\epsilon$. 
Define the angular variable $\phi<\pi/2$ to be the angle between the lines $r=1$ and the line through $e_1$ and positive $z$ axis. Also, for any $0<\gamma<\pi/2$, define $D_1^\gamma$ to be the intersection of $D$ with the sector $\pi/2-\gamma\ge \phi\ge 0$. Denote $D_2^\gamma$  to be the intersection of $D$ with the sector $\pi/2\ge \phi\ge \gamma$ and $D^+$.

\medskip
The following lemma will be key in proving Theorem \ref{lower}.

\begin{lemma}
\label{mainlemma1} There exists a small $\delta>0$ such 
that for all $x:=(\bar{r},\bar{z})\in D_1^\gamma$ with $|x-e_1|<\delta$,
\begin{align}
\label{uzlaw}
u^z(x) = - \frac{4}{\pi}\bar{z}\cdot \int_{Q(\bar{r},\bar{z})} 
\frac{(1-r)z}{((1-r)^2+z^2)^2}w(r,z)\, dr\, dz+ \bar{z} B_1(\bar{r},\bar{z},t)
\end{align}
where $|B_1(\bar{r},\bar{z},t)| \le C(\gamma) \|w_0\|_{L^\infty}$.
Similarly, for all $x\in D_2^\gamma$ with $|x-e_1|<\delta$,
\begin{align}
\label{urlaw}
u^r(x)= \frac{4}{\pi}(1-\bar{r})\cdot \int_{Q(\bar{r},\bar{z})} 
\frac{(1-r)z}{((1-r)^2+z^2)^2}w(r,z)\, dr\, dz+ (1-\bar{r}) B_2(\bar{r},\bar{z},t)
\end{align}
where $|B_2(\bar{r},\bar{z},t)| \le C(\gamma) \|w_0\|_{L^\infty}$.
\end{lemma}

\smallskip
\noindent
{\bf Remark:} In the proof of Theorem \ref{upper}, we needed to carefully keep track of the powers of $r$ and $\bar{r}$. However, since we now are examining dynamics in a neighborhood of $(r,z)=(1,0)$ which is away from the axis, in many cases, powers of $r$ and $\bar{r}$ can safely be controlled by uniform constants.

\Pf Let us prove the expansion for $u^z$ as the one for $u^r$ can be done similarly. For $x=(\bar{r},\bar{z})\in D_1^\gamma$ with $|x-e_1|<\delta$, we have $1-\bar{r} \le (\cot \gamma)\bar{z}$. Define
\begin{align*}
\rho = 10(1+\cot\gamma) \bar{z}
\end{align*}
so we are assured that $x\in B_\rho(e_1)$. Pick $\delta<1/2$ small enough such that $B_\delta(e_1) \subset S_N$ and $\rho<N/2$. The part of the integral for $u^z$ over the ball $B_\rho(e_1)$ will be the part of the remainder term. Using the bounds $|F'(s)|\lesssim 1/s$ and $\frac{r^{3/2}}{\sqrt{\bar{r}}} \lesssim 1$, one can get that this integral is bounded by
\begin{align*}
C \|w_0\|_{L^\infty} \int_{D^+\cap B_\rho(e_1)} \frac{r^{3/2}}{\sqrt{\bar{r}}} \frac{1}{|x-y|} \, dy \le C(\gamma) \|w_0\|_{L^\infty}\bar{z}.
\end{align*}
Next, we will estimate the integral on $S_N\setminus B_\rho(e_1)$ and will remark on $D^+\setminus S_N$ 
later.
On the set $S_N\setminus B_\rho(e_1)$, we can use the Taylor expansion of $F(s)$ at $s=0$ to get 
\begin{align}
\label{GTaylor}
\frac{1}{2\pi}G(r,z,\bar{r},\bar{z}) &= -\frac{1}{4\pi}\sqrt{r\bar{r}} \left[ \log\left(\frac{|x-y|^2}{r\bar{r}}\right) -
\log\left( \frac{|x-y^\ast|^2}{r^\ast\bar{r}}\right)-\log\left(\frac{|\wtilde{x}-y|^2}{r\bar{r}}\right)+\log\left(\frac{|\wtilde{x}-y^\ast|^2}{r^\ast\bar{r}}\right)\right]\\ &\, + O\left(\frac{|x-y^\ast|^2}{r^\ast\bar{r}}\log\left(\frac{|x-y^\ast|^2}{r^\ast\bar{r}}\right)\right)+O\left(\frac{|x-y|^2}{r\bar{r}}\log\left(\frac{|x-y|^2}{r\bar{r}}\right)\right) \nonumber  \\
&= -\frac{1}{4\pi} \sqrt{r\bar{r}} \left[ \log|x-y|^2-\log|x-y^\ast|^2-\log|\wtilde{x}-y|^2+\log|\wtilde{x}-y^\ast|^2\right] \nonumber \\ &\, + O\left(\frac{|x-y^\ast|^2}{r^\ast\bar{r}}\log\left(\frac{|x-y^\ast|^2}{r^\ast\bar{r}}\right)\right)+O\left(\frac{|x-y|^2}{r\bar{r}}\log\left(\frac{|x-y|^2}{r\bar{r}}\right)\right) \nonumber.
\end{align}

\smallskip
Now, we concern ourselves with the four logarithms above as the other terms can be absorbed into the remainder term of \eqref{uzlaw}, which we will remark on later below. We proceed as in \cite{KiselevSverak}. The four logarithms can be written as $-\frac{1}{4\pi}\sqrt{r\bar{r}}$ multiplied with
\begin{align}
\label{logs}
\log \left(1- \frac{2(y-e_1)\cdot(x-e_1)}{|y-e_1|^2} + \frac{|x-e_1|^2}{|y-e_1|^2}\right) - \log\left(1- \frac{2(y^\ast-e_1)\cdot(x-e_1)}{|y^\ast-e_1|^2} + \frac{|x-e_1|^2}{|y^\ast-e_1|^2}\right) \\
-\log \left(1- \frac{2(y-e_1)\cdot(\wtilde{x}-e_1)}{|y-e_1|^2} + \frac{|x-e_1|^2}{|y-e_1|^2}\right)+\log\left(1- \frac{2(y^\ast-e_1)\cdot(\wtilde{x}-e_1)}{|y^\ast-e_1|^2} + \frac{|x-e_1|^2}{|y^\ast-e_1|^2}\right). \nonumber
\end{align} 
We use the following expansion for $\log$ for small $q$
\begin{align*}
\log(1+q)= q-\frac{q^2}{2}+O(q^3).
\end{align*}
On the complement of $B_\rho(e_1)$, $|y-e_1|\ge 10 |x-e_1|$ so we can use this expansion. Then $\eqref{logs}$ becomes
\begin{align*}
-4\frac{z\bar{z}}{|y-e_1|^2}+4\frac{\bar{z} z^\ast}{|y^\ast-e_1|^2}-8\frac{(r-1)(\bar{r}-1)z\bar{z}}{|y-e_1|^4}+
8\frac{(r^\ast-1)(\bar{r}-1)z^\ast\bar{z}}{|y^\ast-e_1|^4}+O\left(\frac{|x-e_1|^3}{|y-e_1|^3}\right).
\end{align*}
Using the identities $\eqref{identities}$ to simplify we get
\begin{align*}
-16\frac{(r-1)(\bar{r}-1)z\bar{z}}{|y-e_1|^4}-8\frac{(\bar{r}-1)z\bar{z}}{|y-e_1|^2} + O\left(\frac{|x-e_1|^3}{|y-e_1|^3}\right).
\end{align*}
Then
\begin{align*}
\frac{1}{4\pi}G(r,z,\bar{r},\bar{z})= \frac{4}{\pi}\sqrt{r\bar{r}}\frac{(r-1)(\bar{r}-1)z\bar{z}}{|y-e_1|^4}+\frac{2}{\pi}\sqrt{r\bar{r}}\frac{(\bar{r}-1)z\bar{z}}{|y-e_1|^2}+\sqrt{r\bar{r}}O\left(\frac{|x-e_1|^3}{|y-e_1|^3}\right).
\end{align*}
Differentiating the above expression with respect to $\bar{r}$ we get
\begin{align}
\label{Gdiff}
\frac{4}{\pi}\sqrt{r\bar{r}} \frac{(r-1)z\bar{z}}{|y-e_1|^4}+ \frac{2}{\pi}\sqrt{r\bar{r}}\frac{z\bar{z}}{|y-e_1|^2} +\frac{2}{\pi}\sqrt{\frac{r}{\bar{r}}} \frac{(r-1)(\bar{r}-1)z\bar{z}}{|y-e_1|^4}+\frac{1}{\pi} \sqrt{\frac{r}{\bar{r}}} \frac{(\bar{r}-1)z\bar{z}}{|y-e_1|^2}+O\left(\frac{|x-e_1|^2}{|y-e_1|^3}\right).
\end{align}
The error term is controlled with
\begin{align*}
|x-e_1|^2 \int_{S_N\setminus B_\rho(e_1)} \frac{1}{|y-e_1|^3}\, dy \lesssim |x-e_1|^2  \int_{\rho}^{1} \frac{1}{t^2}\, dt \lesssim |x-e_1|^2 \rho^{-1} \le C(\gamma)\bar{z}.
\end{align*}
In addition, we can control some of the other terms of \eqref{Gdiff} by
\begin{align*}
\bar{z}\int_{S_N\setminus B_\rho(e_1)} \left(\frac{\sqrt{r\bar{r}}z}{|y-e_1|^2}+\sqrt{\frac{r}{\bar{r}}}\frac{|r-1||\bar{r}-1|z}{|y-e_1|^4}+\sqrt{\frac{r}{\bar{r}}} \frac{|\bar{r}-1|z}{|y-e_1|^2}\right)\, dy\\ \le C\bar{z} \int_{S_N\setminus B_\rho(e_1)}\left(\frac{z}{|y-e_1|^2}+\frac{|r-1||\bar{r}-1|z}{|y-e_1|^4}+\frac{|\bar{r}-1|z}{|y-e_1|^2}\right)\, dy \\
\le C\bar{z} \int_{S_N\setminus B_\rho(e_1)} \left(\frac{1}{|y-e_1|}+1\right)\, dy \le C\bar{z} \int_{\rho}^{1} dt \le C(\gamma) \bar{z}
\end{align*}
Therefore, only the first term of $\eqref{Gdiff}$ will contribute to the main term of \eqref{uzlaw}. Next,
\begin{align*}
\int_{S\setminus B_\rho(e_1)} \frac{(r-1)z}{|y-e_1|^4}\, dy=O(1)+ \int_{Q(\bar{r},\bar{z})} \frac{(r-1)z}{|y-e_1|^4} \, dy
\end{align*}
since we have the bound
\begin{align*}
\left|\int_{Q(\bar{r},\bar{z})\cap B_\rho(e_1)}  \frac{(r-1)z}{|y-e_1|^4} \, dy\right| &\le \left|\int_{\bar{z}}^{C\bar{z}}\int_1^{C\bar{z}} \frac{(r-1)z}{|y-e_1|^4} \, dr\, dz \right|\le C \int_{\bar{z}}^{C\bar{z}}\int_0^{(C\bar{z}-1)^2} \frac{z}{(t+z^2)^2}\, dt\, dz \\
&\le C\int_{\bar{z}}^{C\bar{z}} \frac{1}{z}\, dz \le C.
\end{align*}
Now, $S_N\setminus (Q(\bar{r},\bar{z})\cup B_\rho(e_1))$ is union of two regions, one along the $r$ axis and one near the boundary. These contributions can be controlled using the following bounds
\begin{align*}
\left|\int_{\bar{r}}^1 \int_{\bar{z}}^{N}  \frac{(r-1)z}{|y-e_1|^4}\, dr\, dz\right| \le \int_{\bar{r}}^1 \frac{(1-r)}{(r-1)^2+\bar{z}^2}\, dr\, dz \le C(N)
\end{align*}
\begin{align*}
\left|\int_{N}^{\bar{r}} \int_{0}^{\bar{z}}  \frac{(r-1)z}{|y-e_1|^4}\, dr\, dz\right|\le C(N)
\end{align*}

Thus, the proof of $\eqref{uzlaw}$ will be complete as long as we show that the contribution by integration over $D^+ \setminus S_N$ is controlled up to a constant factor by $\bar{z}$. The expansions of $F$ at $0$ are no longer valid. We can just bound these integrals using Lemma \ref{Fest}. After a computation, we have
\begin{align*}
u^z(\bar{r},\bar{z}) = \int_{D^+} \left(\frac{r}{\bar{r}}\right)^{3/2}(\wtilde{\mathscr{J}}(\bar{r},\bar{z},r,z)-\wtilde{\mathscr{J}}(\bar{r},-\bar{z},r,z)-\wtilde{\mathscr{J}}(\bar{r},\bar{z},r^\ast,z^\ast)+\wtilde{\mathscr{J}}(\bar{r},-\bar{z},r^\ast,z^\ast))w(\bar{r},\bar{z},r,z)\, dr\, dz
\end{align*}
On region $D^+ \setminus S_N$, let us bound $ \left(\frac{r}{\bar{r}}\right)^{3/2}(\wtilde{\mathscr{J}}(\bar{r},\bar{z},r,z)-\wtilde{\mathscr{J}}(\bar{r},-\bar{z},r,z)) = \wtilde{J}(\bar{r},\bar{z},r,z)-\wtilde{J}(\bar{r},-\bar{z},r,z)$. The contribution from the other difference of $\wtilde{\mathscr{J}}$'s will be similar. It suffices to bound
\begin{align*}
\frac{(\bar{r}-r)\sqrt{r}}{\bar{r}^{3/2}}\left(F'\left(\frac{(r-\bar{r})^2+(z-\bar{z})^2}{\bar{r}r}\right) - F'\left(\frac{(r-\bar{r})^2+(z+\bar{z})^2}{\bar{r}r}\right)\right).
\end{align*}
This can be achieved by using Lemma \ref{Fest} and the Mean Value theorem
\begin{align*}
\left|F'\left(\frac{(r-\bar{r})^2+(z-\bar{z})^2}{\bar{r}r}\right) - F'\left(\frac{(r-\bar{r})^2+(z+\bar{z})^2}{\bar{r}r}\right)\right| &\le 2\bar{z} \sup_{\wtilde{z}\in (-\bar{z},\bar{z})} \left|\frac{\partial}{\partial \wtilde{z}} F'\left(\frac{(r-\bar{r})^2+(z+\wtilde{z})^2}{\bar{r}r}\right)\right| \\
& \le C\bar{z} r\bar{r} \frac{1}{(r-\bar{r})^3}.
\end{align*}
Thus, as the other terms satisfy similar bounds, we get
\begin{align*}
\left|\int_{D^+\setminus S_N}\wtilde{J}(\bar{r},\bar{z},r,z)-\wtilde{J}(\bar{r},-\bar{z},r,z)w\, dr\, dz\right| \le C\|w_0\|_\infty\bar{z} \int_{D^+\setminus S_N} \frac{1}{(r-\bar{r})^3}\, dr\, dz \le C(N)\|w_0\|_\infty\bar{z}
\end{align*}

Using a similar argument as directly above, one can show that the terms gotten by formally differentiating the error terms of $\eqref{GTaylor}$
can be similarly controlled by $C\|w_0\|_\infty \bar{z}$.

\medskip
Thus, combining the estimates done above we have that for $x=(\bar{r},\bar{z})\in D_1^\gamma$ with $|x-e_1|<\delta$,
\begin{align*}
u^z(x) &= - \frac{4}{\pi}\frac{\bar{z}}{\sqrt{\bar{r}}}\int_{Q(\bar{r},\bar{z})} 
\frac{r^{3/2} (1-r)z}{((1-r)^2+z^2)^2}w(r,z)\, dr\, dz+ \bar{z} B_1(\bar{r},\bar{z},t) \\
&= - \frac{4}{\pi}\bar{z}\int_{Q(\bar{r},\bar{z})} 
\frac{(1-r)z}{((1-r)^2+z^2)^2}w(r,z)\, dr\, dz+ \bar{z} B_1(\bar{r},\bar{z},t)
\end{align*}
where $|B_1(\bar{r},\bar{z},t)| \le C(\gamma) \|w_0\|_{L^\infty}$. In the last equality above, we can remove 
some factors of $\bar{r}$ and $r$ as the errors they produce can be controlled by the error term $\bar{z} B_1$.

\smallskip
The proof of the lemma is complete $\Box$.

\bigskip
The proof of {\bf Theorem \ref{lower}} now follows from a similar argument as seen in Kiselev and Sverak \cite{KiselevSverak} using Lemma \ref{mainlemma1} in place of their Lemma 3.1. We sketch the details below.

\bigskip
\subsection{Proof of Theorem \ref{lower}}

We start with smooth initial data $w_0$ which is identically $1$ on $D^+$ except on a 
strip of width $\delta$ around $z=0$ where $0<w_0(x)<1$. We assume $w_0$ is odd with respect
to $z$ so $w_0=0$ on $z=0$. Below, we will impose more restrictions on $w_0$. As $u$ is
 incompressible \eqref{divfree}, the distribution function of $w(r,z,t)$ (with respect to the measure
$r\, drdz$) remains the same for all time. This implies the measure of the region where $0<w_0<1$ does not exceed $2\delta$.

\medskip
Then for $|x-e_1|<\delta$ and $x\in D^+$, we can bound the integral term appearing in $\eqref{urlaw}$ and $\eqref{uzlaw}$ by
\begin{align*}
\int_{Q(\bar{r},\bar{z})} 
\frac{(1-r)z}{((1-r)^2+z^2)^2}w(r,z)\, dr\, dz \ge c_1 \int_{c_2\sqrt{\delta}}^{N/2} \int_{\pi/6}^{\pi/3} \frac{1}{\rho}\, d\phi d\rho \ge c_2 \log \delta^{-1}
\end{align*}
Here $c_1$ and  $c_2$ are universal postive constants and we can choose $\delta$ small enough (dependent on $N$) such that the rightmost inequality above holds. The coordinates $(\rho,\phi)$ are polar coordinates centered at $e_1$.  Now, if necessary,
we can choose $\delta$ smaller such that $c_2 \log \delta^{-1} > 100 \cdot C(\gamma)$ where $C(\gamma)$ is the constant from \ref{mainlemma1}.

\medskip
Let $0<z_1'<z_1'' <1$. Define
\begin{align*}
\mathscr{O} (z_1',z_1'') = \left\lbrace (r,z)\in D^+:\, z > -r+1, \, z_1' < z < z_1''\right\rbrace
\end{align*}
along with the quantities
\begin{align*}
\underline{u^z}(z,t) &= \min\{ u^z(r,z,t): \, (r,z)\in D^+, z > -r+1\} \\
\ol{u^z}(z,t) &= \max\{ u^z(r,z,t): \, (r,z)\in D^+, z > -r+1\}.
\end{align*}
From this we can define quantities $a(t)$ and $b(t)$ by
\begin{align*}
\dot{a}(t) &= \ol{u^z}(a,t), \quad a(0)=\epsilon^{10} \\
\dot{b}(t) &= \underline{u^z}(b,t), \quad b(0) = \epsilon.
\end{align*}
where $0<\epsilon<\delta$. Define $\mathscr{O}_t = \mathscr{O}(a(t),b(t))$. Now, we choose
$w_0$ satisfing the same assumptions as above but also specify $w_0=1$ on $\mathscr{O}_0$ with smooth cutoff to zero
 (for example, $\|\nabla w_0 \|_{L^\infty} \sim \epsilon^{-10}$).

\medskip
Now, with these notations in place, one can proceed exactly as in \cite{KiselevSverak}. Using their arguments, one can show that $\mathscr{O}_t$ 
will be non-empty for all $t>0$ and $w(r,z,t)=1$ on $\mathscr{O}_t$. From this, using lemma \ref{mainlemma1}, one can show
$ a(t) \le \epsilon^{C \exp(Ct)}$ for some positive constant $C$. A particle trajectory starting at $z=\epsilon^{10}$ on $\partial D$ near $e_1$ will
never exceed $a(t)$. From this fact, one can arrive at the main estimate of Theorem \ref{lower}. $\Box$

\section{Acknowledgements}

The author would like to thank Prof. Alexander Kiselev for many helpful comments on the numerous drafts of this manuscript.

\begin{appendix}

\bigskip

\section{Estimates for $F$}

Here we will give a rough derivation of the Taylor expansions for the integral $F(s)$
$$
F(s)=\int_0^\pi \frac{\cos\theta \, d\theta}{\sqrt{2(1-\cos\theta)+s}}.
$$
These expansions can also be found in \cite{Sveraknotes}. We will derive them below and that, while elementary, are nonstandard.

\subsection{Estimates at $0$}

We can write $F$ as 
\begin{align*}
F(s)=\int_0^{\pi/2} \frac{1+2\sigma^2}{\sqrt{\sin^2\varphi+\sigma^2}}\, d\varphi -2\int_0^{\pi/2} \sqrt{\sin^2\varphi +\sigma^2}\, d\varphi, \quad \sigma^2=s/4.
\end{align*}
The leading order term above is 
$$
f(\sigma)=\int_0^{\pi/2} \frac{d\varphi}{\sqrt{\sin^2\varphi+ \sigma^2}} = \int_0^{\pi/2} \frac{\cos \varphi\, d\varphi}{\sqrt{\sin^2\varphi+\sigma^2}}+ \int_0^{\pi/2}  \frac{(1-\cos \varphi)\, d\varphi}{\sqrt{\sin^2\varphi+\sigma^2}} :=I+II.
$$
For $I$ we can compute directly and use Taylor series to get
\begin{align*}
I &= \int_0^1\frac{d\varphi}{\sqrt{\varphi^2+\sigma^2}} =\log\frac{1}{\sigma} + \log(1+\sqrt{1+\sigma^2}) \\  &= \log \frac{1}{\sigma} +\log 2+ O(\sigma^2)= \frac{1}{2}\log \frac{1}{s} +2\log 2+ O(s).
\end{align*}

Similarly for $II$ we can get that for $\sigma\to 0^+$
\begin{align*}
II &= \int_0^{\pi/2} \frac{1-\cos\varphi}{\sin\varphi}\, d\varphi+ \int_0^{\pi/2} (1-\cos\varphi)\left(\frac{1}{\sqrt{\sin^2\varphi +\sigma^2}}-\frac{1}{\sin\varphi}\right) d\varphi\\
&= \log 2+ \sigma^2 \int_0^{\pi/2} \frac{1-\cos\varphi}{\sin\varphi\sqrt{\sin^2\varphi+\sigma^2}(\sqrt{\sin^2\varphi+\sigma^2}+\sin\varphi)}d\varphi\\
&= \log 2 + O\left( \sigma^2 \log \frac{1}{\sigma}\right) = \log 2+ O\left( s\log \frac{1}{s}\right)
\end{align*}
where we use that
\begin{align*}
\int_0^{\pi/2} \frac{1-\cos\varphi}{\sin\varphi\sqrt{\sin^2\varphi+\sigma^2}(\sqrt{\sin^2\varphi+\sigma^2}+\sin\varphi)}d\varphi &\le \int_0^{\pi/2} \frac{1-\cos\varphi}{\sin^2\varphi} \frac{d\varphi}{\sqrt{\sin^2\varphi+\sigma^2}} \\
&\le \int_0^{\pi/2} \frac{1-\cos\varphi}{\sin^2\varphi} \frac{d\varphi}{\sqrt{\varphi^2/4+\sigma^2}} \\
&= O\left(\log \frac{1}{\sigma}\right).
\end{align*}
Similarly, we can also have
\begin{align*}
2\int_0^{\pi/2} \sqrt{\sin^2\varphi +\sigma^2}\, d\varphi &= 2+2\sigma^2 \int_0^{\pi/2} \frac{1}{\sqrt{\sin^2\varphi+\sigma^2}+\sin\varphi}\, d\varphi\\
&=2+O\left(\sigma^2 \log \frac{1}{\sigma}\right) = 2+ O\left(s \log \frac{1}{s}\right).
\end{align*}
Putting these expressions together, we get the desired expansion for $F$:
$$
F(s) = \frac{1}{2} \log \frac{1}{s} + \log 8 - 2 + O\left(s \log \frac{1}{s}\right).
$$
Now, consider the derivative
$$
F'(s)= -\frac{1}{2} \int_0^\pi  \frac{\cos \theta \, d\theta}{(2(1-\cos \theta)+s)^{3/2}}.
$$
With $\sigma$ as set above, 
\begin{align}
\label{Fprime}
F'(s) = -\frac{1}{8}\int_0^{\pi/2} \frac{1+2\sigma^2}{(\sin^2\varphi+\sigma^2)^{3/2}}\, d\varphi+ \frac{1}{4}\int_0^{\pi/2} \frac{1}{\sqrt{\sin^2\varphi+\sigma^2}}\, d\varphi.
\end{align}
Doing a similar decomposition as in $F$ above
\begin{align}
\label{Fprime1}
\int_0^{\pi/2} \frac{d\varphi}{(\sin^2\varphi+\sigma^2)^{3/2}}&=\int_0^{\pi/2} \frac{\cos\varphi d\varphi}{(\sin^2\varphi+\sigma^2)^{3/2}}+ \int_0^{\pi/2}\frac{(1-\cos\varphi)d\varphi}{(\sin^2\varphi+\sigma^2)^{3/2}} \\
&= \int_0^{1/\sigma} \frac{dt}{\sigma^2(t^2+1)^{3/2}}+O(\log s) = \frac{4}{s} + O(\log s) \nonumber
\end{align}
where for the second integral on the right hand side above we estimate as we did II above.
The expansion for the second integral in $\eqref{Fprime}$ is done above and putting things together
$$
F'(s)= -\frac{1}{2}\frac{1}{s} + O(\log s), \quad s\to 0^+.
$$

\subsection{Estimates at $\infty$}
Write $F$ as 
$$
F(s) =s^{-1/2} \int_0^\pi \cos \theta \left( \frac{2-2\cos\theta}{s}+1\right)^{-1/2}\, d\theta.
$$
Then by Taylor expansion,
\begin{align*}
F(s) &= s^{-1/2} \int_0^\pi \cos \theta\, d\theta - \frac{1}{2}s^{-3/2}\int_0^\pi 2\cos\theta(1-\cos \theta)\, d\theta + O(s^{-5/2}) \\
&= \frac{\pi}{2} s^{-3/2} + O(s^{-5/2}), \quad s\to \infty.
\end{align*}
The expansions for the derivatives can be derived similarly.

\section{Estimates for $\partial_{\bar{r}} u^z$ and $\partial_{\bar{z}} u^r$ in Theorem \ref{kato}}

 In this section, we sketch the details of the proof of Theorem \ref{kato} for the derivatives 
$\partial_{\bar{r}} u^z$ and $\partial_{\bar{z}} u^r$ under the assumption $\mbox{dist}(x, \partial D) >2\delta$.

\medskip
First, let us do $\partial_{\bar{z}}u^r$.
\begin{align*}
\partial_{\bar{z}}u^r(\bar{r},\bar{z}) &=\left( \int_{B_{\delta \bar{r}}(x)} +\int_{\Omega \cap B_{\delta\bar{r}}^c(x)} +\int_{D\setminus \Omega}\right)\partial_{\bar{z}}K(r,z,\bar{r},\bar{z})w(r,z)\, dr\, dz \\ 
&= I+II+III.
\end{align*}
It suffices to bound the integral $I$ as the other two will be controlled just as in the estimates for $\partial_{\bar{r}} u^r$ earlier.
On $B_{\delta \bar{r}}(x)$, we use the expansions for $F$ and we get
\begin{align*}
\partial_{\bar{z}} \wtilde{K}(r,\bar{r},z,\bar{z}) &= -\frac{\sqrt{r}}{\pi\bar{r}^{3/2}} \left[ F'\left(s\right)+2\frac{(z-\bar{z})^2}{\bar{r}r} F''\left(s\right) \right] = \frac{\sqrt{r}}{2\pi\bar{r}^{3/2}}\left[ \frac{1}{s}-2\frac{(z-\bar{z})^2}{\bar{r}r} \frac{1}{s^2} +O(\log s)\right] \\
&=  \frac{r^{3/2}}{2\pi\bar{r}^{1/2}}\left[ \frac{|x-y|^2-2(z-\bar{z})^2}{|x-y|^4} \right] + \frac{\sqrt{r}}{2\pi\bar{r}^{3/2}} O(\log s)
\end{align*}
where as earlier, we set $\ds s:= \frac{(r-\bar{r})^2+(z-\bar{z})^2}{r\bar{r}}$.
Using a similar argument as in $u_r^r$, we can bound the contribution from the first term
\begin{align*}
\left| \int_{B_{\delta \bar{r}}(x)}\frac{r^{3/2}}{2\pi\bar{r}^{1/2}}\left[ \frac{|x-y|^2-2(z-\bar{z})^2}{|x-y|^4} \right]  w(r,z)\, dr\, dz \right| \le  C(\alpha)\bar{r}^{1+\alpha} \|w_0\|_{L^\infty}.
\end{align*}
We use the fact that the integration of the term in brackets over the ball is $0$.
Furthermore, we can bound the contribution from the error
\begin{align*}
\left| \int_{B_{\delta \bar{r}}(x)}\frac{\sqrt{r}}{2\pi\bar{r}^{3/2}} (\log s) \cdot w(r,z)\, dr\, dz \right| \lesssim \bar{r}(1+\log (\delta^{-1} \bar{r}^{-1})) \|w_0\|_{L^\infty}.
\end{align*}

\bigskip

Now, let us do $u_{\bar{r}}^z$. Again, it will suffice to estimate the integral over $B_{\delta \bar{r}}(x)$.

 Recall 
\begin{align*}
\partial_{\bar{r}} \wtilde{J} =\frac{1}{\pi}\frac{\sqrt{r}}{\bar{r}^{3/2}}\left[F'(s)+
2\frac{(r-\bar{r})^2}{\bar{r}r}F''(s)\right]+\frac{2(r-\bar{r})}{\pi} \frac{\sqrt{r}}{\bar{r}^{5/2}}[F'(s)+sF''(s)]+ \frac{1}{\pi} \frac{r^{3/2}}{\bar{r}^{5/2}}\left[- \frac{3}{8}F(s)
+ sF'(s)+\frac{1}{2}s^2 F''(s)\right].
\end{align*}
The first term can be estimated in the same way as for $u_{\bar{z}}^r$. We can bound the other terms as follows
\begin{align*}
\left| \int_{B_{\delta \bar{r}}(x)} \frac{2(r-\bar{r})}{\pi} \frac{\sqrt{r}}{\bar{r}^{5/2}}[F'(s)+sF''(s)]\, w(r,z)\, dr\, dz \right| & \le  \int_{B_{\delta \bar{r}}(x)}  \left(\frac{r}{\bar{r}}\right)^{3/2} \frac{1}{|x-y|} \, |w(r,z)|\, dr \, dz \\ &\le C \bar{r} \|w_0\|_{L^\infty}
\end{align*}
\begin{align*}
\left| \int_{B_{\delta \bar{r}}(x)} \frac{r^{3/2}}{\bar{r}^{5/2}}\left[- \frac{3}{8}F(s)
+ sF'(s)+\frac{1}{2}s^2 F''(s)\right]\, w(r,z)\, dr\, dz \right| &\le  \int_{B_{\delta \bar{r}}(x)} \frac{r^{3/2}}{\bar{r}^{5/2}} \left| \log(s) \cdot w(r,z)\right|\, dr\, dz \\
& \le C \bar{r} (1+\log(\delta^{-1}\bar{r}^{-1}))\|w_0\|_{L^\infty}.
\end{align*}

\end{appendix}

\end{document}